\magnification 1095

\font\tenBbb=msbm10  \font\sevenBbb=msbm7  \font\fiveBbb=msbm5
\newfam\Bbbfam
\textfont\Bbbfam=\tenBbb \scriptfont\Bbbfam=\sevenBbb
\scriptscriptfont\Bbbfam=\fiveBbb
\def\Bbb{\fam\Bbbfam\tenBbb}

\def\N{{\Bbb N}}
\def\supp{\mathop{\rm supp}\nolimits}
\def\e{\hskip 1.5pt}
\def\r#1#2{\hbox{\noindent\hbox to 1.1cm{[#1]\hfill}
		   \vtop{\advance\hsize by -1.1cm \noindent#2}}}

\newskip\ttglue 
\def\eightpoint{\def\rm{\fam0\eightrm}
  	\textfont0=\eightrm \scriptfont0=\sixrm \scriptscriptfont0=\fiverm
  	\textfont1=\eighti  \scriptfont1=\sixi  \scriptscriptfont1=\fivei
  	\textfont2=\eightsy  \scriptfont2=\sixsy  \scriptscriptfont2=\fivesy
	\textfont3=\tenex  \scriptfont3=\tenex  \scriptscriptfont3=\tenex
	\textfont\itfam=\eightit  \def\it{\fam\itfam\eightit}
	\textfont\slfam=\eightsl  \def\sl{\fam\slfam\eightsl}
	\textfont\ttfam=\eighttt  \def\tt{\fam\ttfam\eighttt}
	\textfont\bffam=\eightbf  \scriptfont\bffam=\sixbf
	\scriptscriptfont\bffam=\fivebf  \def\bf{\fam\bffam\eightbf}
	\tt  \ttglue=.5em plus.25em minus.15em
	\normalbaselineskip=11pt
	\setbox\strutbox=\hbox{\vrule height8pt depth3pt width0pt}
	\let\sc=\sevenrm  \let\big=\eightbig \normalbaselines\rm}

	\font\eightrm=cmr8 \font\sixrm=cmr6 \font\fiverm=cmr5
	\font\eighti=cmmi8  \font\sixi=cmmi6   \font\fivei=cmmi5
	\font\eightsy=cmsy8  \font\sixsy=cmsy6 \font\fivesy=cmsy5
	\font\eightit=cmti8  \font\eightsl=cmsl8  \font\eighttt=cmtt8
	\font\eightbf=cmbx8  \font\sixbf=cmbx6 \font\fivebf=cmbx5
	\def\eightbig#1{{\hbox{$\textfont0=\tenrm\textfont2=\tensy
	\left#1\vbox to7.25pt{}\right.$}}}

\centerline{\bf A remark about distortion,}
\centerline{by Bernard Maurey.}
\bigskip

{\eightpoint  In this note we show that every Banach space $X$
not containing $\ell_1^n$ uniformly and
with unconditional basis
contains an arbitrarily distortable subspace.}
\bigskip

 Let $X$ be a Banach space and $\lambda>1$ a real number.
We say that $X$ is {\it $\lambda$-distortable}
if there exists an equivalent norm $|\e .\e |$ on~$X$ such that
for every infinite dimensional vector subspace $Y$ of~$X$,
$$ \sup\{ |y|/|z|; y,z\in Y, \|y\| = \|z\| = 1\} \ge \lambda.$$
We say that $X$ is {\it arbitrarily distortable} if it is
$\lambda$-distortable for every $\lambda>1$.
The space $X$ is {\it distortable} if it is $\lambda$-distortable
for some $\lambda>1$. It is a classical result of
R.C. James that $\ell_1$ and $c_0$ are not distortable
([J]; see also [LT1], 2.e.3). On the other hand, V. Milman proved
that a non-distortable space $X$ must contain some $\ell_p$
or $c_0$ [M2]. He also formulated there  the general
distortion problem, or rather the ``no distortion conjecture", namely
that there exists no distortable space.
This was disproved by Tsirelson's example [T] of a space containing
no subspace isomorphic to $\ell_p$ or $c_0$.
The distortion problem for $\ell_p$ (i.e. the question whether
$\ell_p$, for $1<p<\infty$, is distortable) remained open
for several more years. It was solved by
Odell and Schlumprecht [OS], who actually proved that $\ell_p$
is arbitrarily distortable for every $p$ such that 
$1< p <\infty$.
\medskip

 A subset $A$ of the unit sphere $S(X)$ of $X$ is 
called {\it asymptotic} if every
infinite dimensional vector subspace $Y$ of $X$ almost intersects $A$,
i.e. $\inf\{\|y - a\|; y\in Y, a\in A\} = 0$
(this notion appears without a name in [M1], in connection with
the distortion conjecture). When $X$ has a basis,
it is easy to see that for $A$ to be asymptotic in $X$,
it is enough to check that $A$ almost intersects every
block subspace $Y$ of $X$ (a {\it block subspace} of $X$ is an
infinite dimensional vector subspace of~$X$
generated  by a sequence of successive blocks from the given basis);
this is because every infinite dimensional subspace $Y$ of $X$
almost contains a block subspace $Z$, in the sense that for every
$\varepsilon>0$, we can find $Z$ such that
$d(z,Y) \le \varepsilon \|z\|$ for every $z$ in $Z$.
\medskip

 Our proof relies heavily on the solution by Odell and
Schlumprecht of the distortion problem for
$\ell_2$, as it appear in [OS]. It is shown there that one can
construct a family $(C_n)$ of asymptotic subsets of the unit
sphere of $\ell_2$, with the property that these sets are
almost pairwise orthogonal.
 Using these sets we shall prove that every Banach space
not containing $\ell_1^n$ uniformly and with unconditional basis
contains an arbitrarily distortable subspace;
this applies for example to
the convexified Tsirelson spaces (but it is still unknown
whether $T$ or $T^*$ -- the original Tsirelson space and its dual--
contains an arbitrarily distortable subspace).
The strategy of our proof is as follows:
first, we use a result of Milman-Tomczak, that gives the existence
of either a subspace that is arbitrarily distortable, or of an
asymptotically-$\ell_p$ or asymptotically-$c_0$ subspace
(see definition below).
With our hypothesis we have
necessarily $1<p<\infty$. The remaining part of the proof
will be to transfer the sets $(C_n)$ from $\ell_2$ to an
asymptotically-$\ell_p$ space $Y$ with unconditional basis
and not containing $\ell_1^n$ uniformly. To this end
we will mimic in $Y$ the construction of the sets $(C_n)$.
\medskip

 It is clear that our problem is an  isomorphic problem, which means
that we can replace the given norm on $X$ by an equivalent norm.
In particular we can
renorm the space so that the equivalent new norm
gets some nicer features
(for example, 1-unconditional basis and uniformly convex
and uniformly smooth norm for a space with unconditional basis
not containing $\ell_1^n$ uniformly, see [LT2].)
\bigskip

 In order to
show that a given Banach space $X$ is
arbitrarily distortable, we need to find
some $\delta >0$ such that for every $\varepsilon>0$,
there exist two asymptotic subsets $A$ and $B$ in $S(X)$
and a $\delta$-norming
set $A^*$ for $A$ in $B(X^*)$ such that the action of $A^*$
on~$B$ is less than~$\varepsilon$, i.e.
$|a^*(b)| \le \varepsilon$ for all $a^*\in A^*$,
$b\in B$ 
(we say that $A^*$ is $\delta$-norming for $A$ if
$\sup\{ |a^*(a)|; a^*\in A^*\} \ge \delta\|a\|$ for every
$a\in A$).
Consider on $X$ the following equivalent norm
$$ |x|_\varepsilon = \|x\| + {1\over\varepsilon}
\sup \{|a^*(x)|; a^*\in A^*\}.$$
For every infinite dimensional subspace $Y$ of $X$, we can find
by definition of asymptotic sets a point $a\in A$ and a
point $b\in B$ such that $a$ and $b$ almost belong to $Y$. We get
$\|a\| = \|b\| = 1$ while $|a|_\varepsilon \ge \delta/\varepsilon$
and $|b|_\varepsilon \le 2$. This shows that the norm $|\e .\e|_\varepsilon$
is a $\displaystyle{\delta\over {2 \varepsilon}}$-distortion 
of the original norm on~$X$. Since by assumption $\varepsilon>0$ 
can be arbitrary, we deduce that $X$ is arbitrarily distortable.
\bigskip

 Let $x_1,\ldots,x_n$ be a sequence of vectors in $X$. We say
that it is $C$-equivalent to the $\ell_p^n$-basis if
there exist two numbers $d,D$
with $d>0$ and $D/d \le C$
such that
$$ d \, (\sum_{i=1}^n |c_i|^p)^{1/p} \le
\|\sum_{i=1}^n c_i x_i\| \le D (\sum_{i=1}^n |c_i|^p)^{1/p}$$
for all scalars $(c_i)_{i=1}^n$.
\medbreak

 Suppose that $X$ has a basis $(e_n)_{n=0}^\infty$. The {\it support} of a
vector $x = \sum_i a_i e_i$ is the subset of $\N$ consisting of all 
integers $i$
such that $a_i \ne 0$, it is denoted by $\supp(x)$. We say that $x$
and $y$ are successive if $\max\supp(x) < \min\supp(y)$. We say that
$x$ is supported after $k$ if $k< \min\supp(x)$.
\medskip

 The first part of these
remarks will make use of a result of Milman-Tomczak.
Following them, we say that a Banach space $Y$ with a basis
is asymptotically-$\ell_p$ (resp: aymptotically-$c_0$)
if there exists a constant $C$
such that for every integer $n\ge 1$, there exists 
an integer $N = N(n)$ such
that every sequence of $n$ successive 
normalized blocks supported after $N$ is $C$-equivalent
to the $\ell_p^n$-basis (resp: $\ell_\infty^n$-basis). 

\medskip

\noindent {\bf Theorem (Milman-Tomczak).} {\sl Let $X$ be an infinite
dimensional Banach space. If $X$ does not contain an arbitrarily
distortable subspace, then $X$ contains an asymptotically-$\ell_p$
or asymptotically-$c_0$ subspace.}
\medskip

 For the convenience of the reader, we shall prove this theorem.
Our proof is slightly shorter than that of [MT].
\smallskip

 Since we may replace $X$ by a subspace and use renorming, we
suppose that $X$ is a Banach space with bimonotone 
basis. We also assume that $X$ does not contain any arbitrarily
distortable subspace, and we will show that $X$ contains an
asymptotically $\ell_p$- or $c_0$-subspace.
\medskip

 If $Y$ is a block subspace of $X$, let $(Y>k)$ denote the
subspace of $Y$ consisting of all vectors in $Y$ 
supported after $k$.
 We shall use a simple stabilization principle:
\smallskip

\noindent{\sl Suppose that to every block subspace $Y$ of $X$
is associated a scalar $\alpha(Y) \ge 0$ such that
$Y_1\subset Y$ implies $\alpha(Y_1) \le \alpha(Y)$ and 
$\alpha(Y>k) = \alpha(Y)$ for every $k\ge 1$. Then there exists 
a block subspace $Z$ of $X$ such that $\alpha(Z_1) =  \alpha(Z)$
whenever $Z_1$ is a block subspace such that $Z_1\subset Z$.}
\medskip

  Let us sketch the proof of the above statement:
let $(Z_k)$ be a decreasing sequence
of block subspaces of $X$ such that $Z_0 = X$ and
$$ \alpha(Z_{k+1}) \le \inf\{\alpha(Y); Y\subset Z_k\} + 2^{-k}$$
for every $k\ge 0$. Let $Z$ be a diagonal subspace
of this sequence $(Z_k)$. It is easy to check that $\alpha(Z')
=\alpha(Z)$ for every block subspace $Z'$ of $Z$.
\medskip

 If $(\alpha_m)_{m=0}^\infty$ is a sequence of such functions, let
$$ \alpha(Y) = \sum_{m=0}^\infty 2^{-m}\, \alpha_m(Y)/\alpha_m(X).$$
If $Z$ stabilizes $\alpha$, then $Z$ stabilizes each $\alpha_m$.
In other words, it is possible to stabilize countably many
functions of this type on the same block subspace $Z$.
\medskip

 For every block subspace $Y$ of $X$,
consider the set $K_Y$ of all $p$s such that $Y$ contains
for every integer $m\ge 1$ a sequence of $m$
successive blocks $2$-equivalent to the $\ell_p^m$-basis.
This set is non empty by Krivine's theorem [K],
and it is closed (easy).
If we consider for each rational $q$ the function
$$\alpha_q(Y) = 1 - \max(1, d(q,K_Y))$$
and apply the preceding discussion we see that
passing to some block subspace $Y$
we can stabilize this set $K_Y$ (in other words,
$K_Z = K_Y$ for any further block subspace $Z$ of $Y$).
Let $p$ be any fixed element of $K=K_Y$, and 
let $q$ be the conjugate exponent.
For every block subspace $Y$ of $X$
and every integer $n\ge 1$ let $\tilde F_Y(n)$ be the 
smallest constant $C$ such that
$$ \|y_1+\cdots+y_n\| \le C (\sum_{i=1}^n \|y_i\|^p)^{1/p}$$
for all sequences
$y_1,\ldots,y_n$ of successive elements in $B(Y)$;
let $\tilde G_Y(n)$ be the smallest constant $C$ such that
$$ \sup\{ (x^*_1+\cdots+x^*_n).y \,;\, y\in B(Y)\} \le C
(\sum_{i=1}^n \|x^*_i\|^q)^{1/q}$$
for all sequences $x^*_1,\ldots,
x^*_n$ of successive elements in $B(X^*)$, i.e. successive
blocks of the biorthogonal functionals.
 Let now
$$ F_Y(n) = \lim_k \tilde F_{(Y>k)}(n), 
\ G_Y(n) = \lim_k \tilde G_{(Y>k)}(n).$$
We can stabilize $F_Y$ and $G_Y$ by passing to some 
further block subspace $Y$. It is clear that $F_Y$ and $G_Y$ are
non-decreasing functions of $n$,
with $F_Y(1) = G_Y(1) = 1$.
\smallskip


 From now on we fix a stabilizing block subspace $Y$; we
set $F=F_Y$, $G=G_Y$. Using 
the stabilization assumption we
know that every block subspace of $Y$ contains
for every integer $m \ge 1$ a sequence of $m$ successive blocks
that is $2$-equivalent to the $\ell_p^m$-basis.
We have assumed that~$X$
does not contain any arbitrarily distortable subspace, hence
$Y$ is not arbitrarily distortable, therefore
there exists $D \ge 1$ such that for every equivalent norm on $Y$,
there exists an infinite dimensional subspace $Z$ of $Y$
on which the new norm is $D$-equivalent to the initial norm.
We will prove now that $Y$ is asymptotically-$\ell_p$.
This is clearly equivalent to showing that $\sup_m F(m)$
and $\sup_m G(m)$ are both finite.
\medskip

 We first show that
$\sup_m G(m) \le 8D$. Working by contradiction, let us assume
that $G(m) > 8D$ for some integer $m$.
Since $G(m) > 1$, it implies that $q>1$, hence $p<\infty$.
We consider the subset $A^*$ of $B(Y^*)$ consisting
of all functionals on $Y$ that can be represented as
$(x^*_1+\cdots+x^*_m)/G(m)$
with successive $x^*_j$s in $X^*$
such that $\sum_{i=1}^m \|x^*_i\|^q \le 1$.
Let~$A$ be the set of norm one vectors in $Y$ 
that are $1/4$-normed by
some element in $A^*$, i.e. the set of all 
vectors $a\in S(Y)$ such that
$\sup\{ |a^*(a)|; a^*\in A^*\} \ge \|a\|/4$.
For every block subspace $Z$ of $Y$,
we have by stabilization $G_Z(m) = G(m)$,
hence there exist $z\in S(Z)$ and $x^*_1,\ldots,x^*_m$ 
successive in $B(X^*)$ such that
$\sum_{i=1}^m \|x^*_i\|^q \le 1$,
$$\|x^*_1+\cdots+x^*_m\|_{Y^*} \le 2 G(m)
\hbox{\ \ and\ \ }
(x^*_1+\cdots+x^*_m).z \ge  G(m)/2;$$
this shows that $z\in A$, hence $A$ is
asymptotic in $Y$. Let $N>>m$ and let~$B$ be the set of
$\ell_p^N$-vectors with constant $2$ in $Y$
(i.e. normalized vectors of the form $(y_1+\cdots+y_N)/N^{1/p}$,
where $y_1,\ldots,y_N$ are successive and
$2$-equivalent to the $\ell_p^N$-basis.)
It follows from our assumptions that $B$ is asymptotic
in $Y$, and by an easy argument
--used by Schlumprecht (see [S] or [GM], Lemma~4) in the context
of Schlumprecht's space-- the action of
$A^*$ on $B$ is bounded by $2/G(m)$.
Indeed, if $b=(y_1+\cdots+y_N)/N^{1/p}\in B$ and
$a^* = (x_1^*+\cdots+x_m^*)/G(m)\in A^*$, then up to a small
perturbation of $b$ (bounded by
$2(m/N)^{1/p}$, hence small since $p<\infty$)
the action of the $x_j^*$s cuts~$b$ into
$m$ pieces corresponding to successive
subsets $A_1,\ldots,A_m$ of $\{1,\ldots,N\}$. Then
for each $j=1,\ldots,m$
$$ |x^*_j.(\sum_{i=1}^N y_i)| \simeq |x^*_j.(\sum_{i\in A_j} y_i)| \le
\|x^*_j\| \, \|\sum_{i\in A_j} y_i\| \le 
2 \|x^*_j\| \, |A_j|^{1/p}$$
and we apply H\"older's inequality
$\sum_{j=1}^m  \|x^*_j\| |A_j|^{1/p} \le N^{1/p} $ to
get $|a^*(b)| \le 2/G(m)$.
If we define a new norm on $Y$ which is essentially
given by the supremum on $A^*$, namely
$$ |y| = \sup\{ |a^*(y)|; a^*\in A^*\} + \tau \|y\|$$
(for some small $\tau>0$) then
any block subspace of $Y$ contains elements
in~$B$, whose new norm will not exceed $\tau + 2/G(m)$, while the
new norm of elements in $A$ is at least $1/4$.
This shows that the equivalent new norm is 
distorted by more than $G(m)/(8 + 4\tau G(m)) > D$
(for $\tau$ small enough),
contradicting our assumption
of bounded distortion for $Y$.
\bigskip

 We know now that $C = \sup_m G(m) \le 8D$;
we are going to show that
$\sup_m F(m) \le 4CD$. If not, let $m$ be such that
$F(m) > 4CD$. This implies that
$p>1$ and $q<\infty$.  Let $N>>m$ and let $k$ be such that
$\tilde G_{(Y>k)}(N) \le 2 G(N) \le 2C$.
Let $B$ be the set of all vectors in $B(Y)$
of the form $(y_1+\cdots+y_m)/F(m)$ where the $y_i$s are
successive vectors in $Y$ and $\sum_{i=1}^m \|y_i\|^p \le 1$;
Let~$A$ be the set of $\ell_p^N$-vectors
with~constant $2$ in $(Y>k)$.
If $(x_1+\cdots+ x_N)/N^{1/p} \in A$ 
and $x^*_j\in X^*$ is a norming functional for $x_j$
when $j=1,\ldots,N$, $x^*_j$ supported on the smallest interval
supporting $x_j$, the action of $(x^*_1+\cdots+x^*_N)/N^{1/q}$
on $(x_1+\cdots+x_N)/N^{1/p}$ is equal to
$$ {1\over N} \sum_{i=1}^N \|x_i\| \ge 1/2.$$
We take as $A^*$ the set
of all functionals $(x^*_1+\cdots+x_N^*)/N^{1/q}$,
with $x^*_1,\ldots,x^*_N$ successive in $B(X^*)$, and
supported after $k$.
With the same argument  as before
we check that the action
of $A^*$ on $B$ is bounded by $2C/F(m)$.
Indeed, if $b=(y_1+\cdots+y_m)/F(m)\in B$ and
$a^* = (x_1^*+\cdots+x_N^*)/N^{1/q}\in A^*$, then up to a small
perturbation (using $q<\infty$)
the action of the $y_j$s cuts $a^*$ into
$m$ pieces corresponding to successive
subsets $A_1,\ldots,A_m$ of $\{1,\ldots,N\}$. Then
for each $j=1,\ldots,m$, we obtain
using $\tilde G_{(Y>k)}(N) \le 2 C$ that
$$ |(\sum_{i=1}^N x^*_i).y_j| \simeq |(\sum_{i\in A_j} x^*_i).y_j| \le 
2 C |A_j|^{1/q} \, \|y_j\|$$
and we apply H\"older's inequality
$\sum_{j=1}^m \|y_j\| \, |A_j|^{1/q} \le N^{1/q}$ to
get $|a^*(b)| \le 2C/F(m)$.
By definition of $F(m)$ and stabilization,
every block subspace of $Y$ contains
elements of $B$ that have norms tending to $1$.
If we define as before a new norm on $Y$ 
by taking the supremum on $A^*$,
the new norm of elements in $B$
will not exceed $2 C/F(m)$, while the
new norm of elements in $A$ is at least $1/2$.
This shows that the equivalent norm is 
distorted by more than $F(m)/(4C) > D$,
contradicting again our initial assumption.
At this point we have proved that $Y$ is asymptotically-$\ell_p$.
This ends the proof of the theorem.
\bigbreak

 Here is the main result of this Note:
\smallskip

\noindent{\bf Remark.} {\sl Let $X$ be a Banach space with unconditional
basis and not containing $\ell_1^n$ uniformly. Then $X$
contains an arbitrarily distortable subspace.}
\medskip

 It is clear from the proof of the Theorem that
if a Banach space $X$ with a basis contains no 
arbitrarily distortable subspace, then we find an
asymptotically-$\ell_p$ subspace $Y$ as a
block subspace of $X$. If we start with a space $X$
with unconditional basis, then
$Y$ has an unconditional basis,
and if $X$ does not contain $\ell_1^n$ uniformly,
then $Y$ of course does not contain $\ell_1^n$ uniformly.
All we need is to show that $Y$ is arbitrarily distortable.
Actually, we are going to prove slightly more.
\medbreak

 A Banach space $X$ is {\it sequentially arbitrarily distortable}
([OS])
if there exists a sequence $|\e .\e|_i$ of equivalent norms on $X$
and a sequence $(\varepsilon_i)$ decreasing to $0$ such that

 -- $ |\e .\e |_i \le \|.\|$ \hbox{\ for each\ } $i$ and

 -- for every subspace $Y$ of $X$ and for every $i_0$ there exists
$y\in S(Y, |\e .\e|_{i_0})$ such that $|y|_i \le \varepsilon_{\min(i, i_0)}$
for $i\ne i_0$.

\medskip

\noindent {\bf Remark 2.} {\sl Let $Y$ be an asymptotically
$\ell_p$-space with unconditional basis and 
not containing $\ell_1^n$ uniformly.
Then $Y$ is sequentially arbitrarily distortable.}
\medskip

 This remark applies for example to
the convexified Tsirelson spaces (see [CS]).
\smallskip

 We can first renorm $Y$ in such a way that the basis is
1-unconditional and the norm uniformly convex and uniformly smooth
(see [LT2], section 1f).
Since $Y$ is asymptotically-$\ell_p$, there exists $C$ such that
for every integer $n\ge 1$,
every sequence of $n$ successive 
normalized blocks ``far enough"
(depending upon $n$) is $C$-equivalent to the $\ell_p^n$-basis,
and the same in $Y^*$ with the $\ell_q^n$-basis. Precisely, to every
integer $n\ge 1$ we can associate an integer $P(n)$ such that
$$ C^{-p}(\sum_{i=1}^n \|y_i\|^p) \le \|\sum_{i=1}^n y_i\|^p \le
C^p (\sum_{i=1}^n \|y_i\|^p)$$
$$ C^{-q}(\sum_{i=1}^n \|y^*_i\|^q) \le \|\sum_{i=1}^n y^*_i\|^q \le
C^q (\sum_{i=1}^n \|y^*_i\|^q)$$
for all sequences of successive elements $y_1,\ldots,y_n$ in
$Y$ or $y^*_1,\ldots,y^*_n$ in $Y^*$ that are supported after $P(n)$.
 We shall refer
to this property as property $A_{n}$ in what follows.
It is clear that $1< p < \infty$, since $Y$ 
does not contain $\ell_1^n$ uniformly.
\smallskip

 The construction in [OS] uses properties of the space
$S$ constructed by Schlumprecht in [S].
 It is perhaps useful for the reader to recall a minimal amount
of information about the space~$S$.
Let $\varphi(x) = \log_2(x+1)$ for $x\ge 0$.
Let $B$ be the smallest
subset of $B(c_0)$ containing $\pm e_i$ for every basis vector
of $c_0$, $i\in \N$, and such that for every
integer $n\ge 1$, when $v_1,\ldots,v_n$ are
successive elements in $B$, then
$(v_1+\cdots+v_n)/\varphi(n)$ belongs to $B$. Define a norm
on $\ell_1$ by
$$ \|u\|_S = \sup \{ |v(u)|; v\in B\}$$
and let $S$ be the completion of $\ell_1$ for this norm.
We recall two important facts about $S$:
the canonical basis is 1-unconditional in $S$; using the
notation introduced in the proof of the Theorem, we have that
the Krivine set $K_S$ for Schlumprecht's space
reduces to $\{1\}$.
\medskip

 Since $Y$ has a 1-unconditional basis, we can consider
$Y$ and $Y^*$ as  lattices of functions on $\N$.
We shall do the same for Schlumprecht's space $S$ and
its dual $S^*$. When $x$ is a vector in $Y$ or $S$ or $S^*$
and $E$ a subset of $\Bbb N$, we denote by $E x$ the projection
of $x$ on the subset $E$. If $u$ is a vector in $S$ and $v\in S^*$,
we consider the pointwise product $uv$: this is the vector of
$\ell_1$ such that $(uv)_i = u_i v_i$ for every $i\ge 0$;
we also
use sometimes the notation $u.v$ instead of $v(u)=
\sum_{i\ge 0} u_i v_i$.
We say that
a vector $u$ is non-negative if $u_i\ge 0$ for every $i\ge 0$.
\medskip

 We need some notation about the space $S$; we say that $x\in S$ is a
$\ell_1^n$-average if $\|x\| = 1$ and $x = n^{-1}
(x_1+\cdots+x_n)$, where $x_1,\ldots,x_n$ is
a sequence of successive vectors, $2$-equivalent to the
$\ell_1^n$-basis. Next, a sequence $u_1,\ldots,u_k$ of
successive vectors in $S$ is called a {\it rapidly increasing
sequence} of $\ell_1$-averages 
(in short: RIS) of length $k$ and sizes $n_1,\ldots,n_k$
if every $u_i$, $i=1,\ldots,k$
is a $\ell_1^{n_i}$-average and the sequence $(n_i)_{i=1}^k$ satisfies the
following lacunarity condition:
$$\varphi(n_1/4) \ge 36 \, k^2 \hbox{\  and\ }
\varphi(n_j)^{1/2} \ge 2 |\supp(u_{j-1})|
\hbox{\ for\ } j=2,\ldots,k.$$
\medskip
 
 For every integer $k\ge 1$, we~denote by $A_k$ the subset of $S$
corresponding to vectors of the form
$u = {{\varphi(k)}\over k} (u_1+\cdots+u_k)$
where $u_1,\ldots,u_k$ is a RIS of length $k$. We call $A^*_k$
the set of all elements in $B(S^*)$ of
the form $v = (v_1+\cdots+v_k)/\varphi(k)$,
where $v_1,\ldots,v_k$ are successive elements in $B(S^*)$.
It is proved in [S] or [GM] that 
\medskip

\noindent{\bf Lemma 1.} {\sl There exists a sequence
$(\varepsilon_k)$ tending to $0$ and a sequence $(L(k))$
of integers such that when
$x^*\in A^*_k$ and $x\in A_l$, or when
$x\in A_k$, $x^*\in A^*_l$, with $l\ge L(k)$,
then $|x^*(x)| \le \varepsilon_{k}$.} 
\medskip

 The fundamental building blocks 
for the distortion of $\ell_2$ in [OS] are the
subsets $D_m$ of the unit sphere
$S(\ell_1)$ of $\ell_1$ consisting of
pointwise products $u v$ where $u$ is an $\ell_1^m$-average
in Schlumprecht's space $S$ and $v\in 2 B(S^*)$
is almost norming $u$. It is proved in [OS]
that for every integer $m \ge 1$ and every $\varepsilon >0$,
there exists $K = K(m,\varepsilon)$ such that given
$h_1,\ldots,h_K$ successive in $S(\ell_1)$, there exists a subset
$J$ of $\{1,\ldots,K\}$, an $\ell_1^m$-average $u\in S$
and $v\in 2 B(S^*)$ such that
$$ \|uv\|_1 = 1, \ \ 
\| uv - {1\over |J|} \sum_{j\in J} h_j\|_1 < \varepsilon.$$
(See [OS], Lemma~3.5 and sublemma~3.6).
Suppose that $Z$ is a block subspace of $Y$. Let
$K = K(m,1/m)$ and let $z_1,\ldots,z_K$ be successive elements in
$S(Z)$, supported after $P(K)$, so that $A_{K}$ applies. Let
$z^*_1,\ldots,z^*_K$ be the norming functionals for these
vectors. Since we assumed $Y$ smooth, each $z_j$ admits
a unique norming functional $z^*_j$;
since the basis of $Y$ is 1-unconditional,
we have $\supp(z^*_j) = \supp(z_j)$
and the pointwise product $z_j z^*_j$ is non-negative.
Apply the result of [OS] to the sequence
$h_j = z_j z^*_j$, $j=1,\ldots,K$. 
We find a subset $J$ of $\{1,\ldots,K\}$,
an $\ell_1^m$-average $u\in S$ and $v\in 2 B(S^*)$ such that
$$ \|u v\|_1 = 1,
\ \| uv - {1\over |J|} \sum_{j\in J} z_j z^*_j \|_1 < {1\over m}.$$
Since the basis in $S$ is also 1-unconditional,
and since we observed that $h_j$ is non-negative,
we may assume that
$u$ and $v$ are non-negative vectors.
Consider now
$$ z = |J|^{-1/p} \sum_{j\in J} z_j, \ \ \
z^* = |J|^{-1/q} \sum_{j\in J} z^*_j.$$
Then $z\in Z$, 
and by $A_{K}$ we get $C^{-1}\le \|z\| \le C$, $C^{-1}\le \|z^*\| \le C$.
Also $z z^* = |J|^{-1} \sum_{j\in J} z_j z^*_j$ therefore
$\| uv - z z^*\|_1 < 1/m$.
\smallbreak

 Let $\Delta_m$ be the subset of $Y\times Y^*$
consisting of all couples $(z,z^*)$
such that  $z = N^{-1/p} \sum_{j=1}^N z_j$,
$z^* = N^{-1/q} \sum_{j=1}^N z^*_j$
for some integer $N\ge 1$,
where $z_1,\ldots,z_N$ are successive in $S(Y)$ and supported
after $P(N)$,
$z^*_j$ is the norming functional for $z_j$,
and there exists an $\ell_1^m$-average $u\in S$
and $v\in 2 B(S^*)$ such that $\|u v\|_1 =1$, $\|uv - z z^*\|_1 < 1/m$,
$u$ and $v$ non-negative.
Notice that $z z^*$ is non-negative and that
$\|z z^*\|_1 = 1$.

 Let $Q$ and $Q^*$ be the projections from $Y\times Y^*$
onto $Y$ and $Y^*$ respectively.  It follows from the
preceding discussion that
for every integer $m\ge 1$, every block subspace
of $Y$ contains an element $z\in Q\Delta_m$.
In other words, this set $Q \Delta_m$ is 
a non-normalized asymptotic set in $Y$.
\smallskip

 Next we recall the construction of the sets $(C_k)$ in [OS]:
actually we rather use subsets $(B_k)$ of $S(\ell_1)$;
each $C_k$ consists of elements whose absolute values
are square root of elements in $B_k$.  The sets $(B_k)$
are built from averages between very different sets
$D_{m_i}$, using the lacunarity conditions for RIS.
More precisely, an element $x$ in $B_k$ as the form
$x = {1\over k}(u_1 v_1 +\cdots +u_k v_k)$, where $u_i v_i\in D_{m_i}$
for $i=1,\ldots,k$ are successive, $u_1,\ldots,u_k$ is a RIS of 
length $k$ and sizes $m_1,\ldots,m_k$ in $S$, and $v_i\in 2 B(S^*)$,
for $i=1,\ldots,k$. This vector $x$ has a representation
as pointwise product
$$ x = \Bigl({{\varphi(k)}\over k} (u_1+\cdots+u_k) \Bigr)
\Bigl((v_1+\cdots+v_k)/\varphi(k) \Bigr)$$
hence $B_k$ appears as a subset of
$ S(\ell_1) \cap 2 A_k A^*_k$.  The sets $(B_k)$
have mutual distances almost $2$ in $\ell_1$, when the indices
are very different; in other words, any element
$x\in B_k$ is almost disjoint from any element $y\in B_l$
when $l \ge L(k)$.
This follows easily from the fact that the
scalar product between $\sqrt{|x|}$ and~$\sqrt{|y|}$ is almost~$0$,
by Lemma 1 (write $|\sqrt{|x|} - \sqrt{|y|}| \le \sqrt{|x - y|}$ and compute
the $\ell_2$-norms), and it will be used at the end of the proof.
\medbreak

 Corresponding to the sets $(B_k)$ we construct now sets 
of couples $(\Gamma_k)$
in $Y\times Y^*$. A couple $(y,y^*)$ belongs to 
$\Gamma_k$ if there exists a sequence of successive couples
$(z_1,z^*_1),\ldots,(z_k,z^*_k)$ supported after $P(k)$
such that $(z_i,z^*_i)\in \Delta_{m_i}$,
and there exists a RIS
$u_1,\ldots,u_k$ of length $k$ with sizes $m_1,\ldots,m_k$,
a sequence $v_1,\ldots,v_k$ in $2 B(S^*)$,
$u_i$ and $v_i$ non-negative, such that
$$ y = k^{-1/p} \sum_{i=1}^k z_i,\ y^* = k^{-1/q} \sum_{i=1}^k z^*_i,$$
$$ \|u_i v_i\|_1 = 1,
\ \|u_i v_i - z_i z^*_i\|_1 < {1\over m_i},\ \forall i=1,\ldots,k.$$
Observe again that $y y^*$ is non-negative and $\|y y^*\|_1 = 1$.
 It follows from
$A_{k}$ that $C^{-2} \le \|y\| \le C^2$, $C^{-2} \le \|y^*\| \le C^2$.
It is clear that for every integer $k\ge 1$, the set
$Q \Gamma_k$ is a non-normalized asymptotic set in $Y$.
\medskip



\noindent{\bf Lemma 2.} {\sl There exists a function
$\alpha(\varepsilon)$ with $\lim_{\varepsilon\rightarrow 0}
\alpha(\varepsilon)  = 0$
such that for every integer $m\ge 1$, if $(y,y^*)\in\Gamma_m$
and $\|E yy^*\|_1 <\varepsilon$ for some subset $E\subset\N$,
then $\|E y\| < \alpha(\varepsilon)\|$ and
$\|E y^*\| < \alpha(\varepsilon)$.}
\medskip

\noindent Proof:  since $Y$ is uniformly convex and uniformly smooth, it 
is clear that there exists a function $\alpha_1(\varepsilon)$
with $\lim_{\varepsilon\rightarrow 0} \alpha_1(\varepsilon) = 0$ such that,
if $y^*\in S(Y^*)$ norms
$y\in S(Y)$, if $E$ is any subset of~$\N$
and if $\|E y y^*\|_1 < \varepsilon$,
then $\|Ey\|, \|Ey^*\| < \alpha_1(\varepsilon)$. Indeed,
if $F$ denotes ${\Bbb N}\setminus E$, we have
$$\|(Fy)y^*\|_1 = \|yy^*\|_1 - \|Eyy^*\|_1 =
1 - \|Eyy^*\|_1 > 1 -\varepsilon,$$ 
hence $\|Fy\| > 1-\varepsilon$;
now $\|Fy \pm Ey\| = 1$, and it follows
by uniform convexity that
$\|Fy\| \le 1 - \delta(\|Ey\|)$, where $\delta$ is the modulus of
convexity of $Y$.
We finally obtain $\delta(\|Ey\|) <\varepsilon$.
The same reasoning applies in the dual.
This implies our claim, taking for $\alpha_1$ the
maximum of the inverse
function of $\delta$ and the inverse function of
$\delta^*$, the modulus of convexity of~$Y^*$.
\smallskip

 Suppose that $(z,z^*)$ is an element of $\Delta_m$,
with $z = N^{-1/p} \sum_{j=1}^N z_j$ and
$z^* = N^{-1/q} \sum_{j=1}^N z^*_j$,
with $z^*_i$
norming $z_i\in S(Y)$ for $i=1,\ldots,N$,
$z_i z^*_i$ successive and supported after $P(N)$, so that
$A_{N}$ applies.
Then $z z^* = N^{-1} \sum_{i=1}^N z_iz^*_i$. Suppose that
$\|E z z^*\|_1 < \varepsilon^2$ for some subset $E$ of~$\N$.
If we consider $J = \{i; \|E z_i z^*_i\|_1 \ge \varepsilon\}$,
we see that $|J| < \varepsilon N$.
For $i\notin J$ we have $\|Ez_i z^*_i\|_1 < \varepsilon$,
hence using the preliminary remark
we deduce that $\| E z_i\| < \alpha_1(\varepsilon)$ for $i\notin J$.
We can then estimate $\|E z\|$, using $A_{N}$
$$ \|Ez\| < C(\alpha_1(\varepsilon) +  (|J|/N)^{1/p}) < 
C(\alpha_1(\varepsilon) + \varepsilon^{1/p}).$$
We can work
similarly in the dual for $\|E z^*\|$.
We deduce that
$$ \|E z z^*\|_1 < \varepsilon^2 \hbox{\ implies\ }
\|E z\| < \alpha_2(\varepsilon),
\ \|E z^*\| < \alpha_2(\varepsilon),$$
where $\alpha_2(\varepsilon) = C(\alpha_1(\varepsilon)+
\varepsilon^{\min(1/p,1/q)})$.
\medskip

 Let $(y, y^*)$ belong to some
$\Gamma_m$, with $y =m^{-1/p} \sum_{i=1}^m y_i$,
$y^* = m^{-1/q} \sum_{i=1}^m y^*_i$,
where the products $(y_i y^*_i)$s are successive,
supported after $P(m)$, and
each $(y_i, y^*_i)$ belongs to some $\Delta _{n_i}$.
By the preceding discussion we have that
$\|E y_i y^*_i\|_1 < \varepsilon^2$ implies $\|E y_i\|,
\|E y^*_i\| < \alpha_2(\varepsilon)$.
Using the same argument and losing an other
factor~$\varepsilon$ we get that $\|E y y^*\|_1 < \varepsilon^3$
implies $\|E y\|, \|E y^*\| < \alpha_3(\varepsilon)$.
This ends the proof of Lemma~2: we simply set
$\alpha(\varepsilon) = \alpha_3(\varepsilon^{1/3})$.
\bigbreak

 We are now in a position to finish the proof of Remark~2.
 Let again $(y, y^*)$ belong to some
$\Gamma_m$, with $y =m^{-1/p} \sum_{i=1}^m y_i$,
$y^* = m^{-1/q} \sum_{i=1}^m y^*_i$,
where the products $(y_i y^*_i)$s are successive,
supported after $P(m)$,
each $(y_i, y^*_i)$ belongs to some $\Delta _{n_i}$,
and
there exist for $i=1,\ldots,m$ vectors
$u_i\in S$, $v_i\in 2 B(S^*)$, $u_i$ and $v_i$ non-negative
such that $u_1,\ldots,u_m$ is a RIS of length $m$
and sizes $n_1,\ldots,n_m$ in $S$,
and finally $\|u_i v_i\|_1 = 1$,
$\|u_i v_i - y_i y^*_i\|_1 < 1/n_i$.
 Let
$$ u ={{\varphi(m)}\over m} \sum_{i=1}^m u_i ;
\ v = \sum_{i=1}^m v_i/\varphi(m) ;$$
 We have
$$ u v = {1\over m} \sum_{i=1}^m u_i v_i; \ \
y y^* = {1\over m} \sum_{i=1}^m y_i y^*_i$$
$$ \|u v\|_1 = 1,\ \| uv - yy^*\|_1 < 1/n_1 < 2^{-36 m^2}.$$
(We used the lacunarity condition for RIS).
We observe that $u\in A_m$ and $v\in 2 A^*_m$.
\smallskip

  In the same way let
$(\tilde y, \tilde y^*)$ belong to $\Gamma_l$,
with $l> L(m)$.
In a similar way $\|\tilde u \tilde v - \tilde y \tilde y^*\|_1<
2^{-36 l^2}$,
for some $\tilde u \in A_l$ and $\tilde v\in 2 A^*_l$. We assumed
$u,v, \tilde u, \tilde v$ non-negative.
If $l \ge L(m)$, we know by Lemma~1
that $\tilde v.u$ and $v.\tilde u$ are smaller than
$2 \varepsilon_m$,
hence $\sqrt{u v}$ is almost orthogonal to
$\sqrt{\tilde u \tilde v}$. Indeed by the Cauchy-Schwarz
inequality
$$ \sqrt{uv} . \sqrt{\tilde u \tilde v}
= \sqrt{u \tilde v}. \sqrt{\tilde u v} \le
\sqrt{u.\tilde v} \, \sqrt{\tilde u.v}
\le 2 \varepsilon_m.$$
 therefore $u v$ and $\tilde u \tilde v$ have distance almost $2$
in $\ell_1$, because letting $t = uv$, $t' = \tilde u \tilde v$,
$$ \|t-t'\|_1  = \|\sqrt{|t-t'|} \, \|_2^2 \ge \|\sqrt t - \sqrt {t'}\|_2^2
= 2 - 2\sqrt t.\sqrt {t'} \ge 2(1- 2 \varepsilon_m).$$
Let us choose $m$ such that $\varepsilon_m < \varepsilon/4$
and $2^{-36 m^2} < \varepsilon/4$. The preceding estimates yield
$$\|yy^* - \tilde y \tilde y^*\|_1 > 2(1- \varepsilon).$$
Now $f = yy^*$ and $\tilde f = \tilde y \tilde y^*$ are non-negative, thus
$$\|f- \tilde f\|_1 = \|f\|_1 + \|\tilde f\|_1 - 2 \|\min(f,\tilde f)\|_1
= 2(1 - \|\min(f,\tilde f)\|_1),$$
hence $\|\min(f, \tilde f)\|_1 < \varepsilon$.
Let $E=\{i\in\N; f_i > \tilde f_i\}$ and
$F =\{i\in\N; f_i \le \tilde f_i\}$.
 We have a partition of $\N$ into two sets $E$
and $F$ such that
$$ \|F y y^*\|_1 <\varepsilon \hbox{\ and\ }
\| E \tilde y \tilde y^*\|_1 < \varepsilon.$$
Using Lemma~2
we deduce that $\| E \tilde y\| < \alpha(\varepsilon)$
and $\|Fy^*\| < \alpha(\varepsilon)$.
We also know that $\|\tilde y\|, \|y^*\| \le C^2$.
Finally 
$$|y^*. \tilde y| = |Fy^*.\tilde y + y^*.E\tilde y| \le
2 C^2 \alpha(\varepsilon).$$
\medskip

 We see that $|\tilde y^*.y|$ is small for the same
reasons; finally, for $l\ge L(m)$, the actions
of $Q^* \Gamma_l$ on $Q \Gamma_m$ and of $Q^* \Gamma_m$ on
$Q \Gamma_l$ are small, depending upon $m$.
This allows to show that $Y$ is sequentially arbitrarily distortable.
To this end let us select a sequence $(m_k)$ of integers
such that
$$ m_{k+1} \ge L(m_k), \ \ \varepsilon_{m_k} < 2^{-k-2},$$
and define for every integer $k\ge 1$ an equivalent norm on $Y$
by the formula
$$ |y|_k = C^{-2} \sup\{ |z^*(y)|; z^*\in Q^*\Gamma_{m_k}\}.$$
First, observe that $\|y\|_k \le \|y\|$ since we know
that $\|z^*\| \le C^2$ when $z^*\in Q^*\Gamma_{m_k}$.
Suppose that $l> k$. By the above proof, we know that the
action of $Q^* \Gamma_{m_k}$ on $Q \Gamma_{m_l}$ and the action
of $Q^* \Gamma_{m_l}$ on $Q \Gamma_{m_k}$ are less than
$ 2C^2 \alpha(2^{-k})$. In other words,
$\|y\|_k \le 2C^2 \alpha(2^{-k})$ when
$y\in Q \Gamma_{m_l}$ and $\|z\|_l \le 2C^2 \alpha(2^{-k})$ when
$z\in Q\Gamma_{m_k}$. On the other hand, $z^*(z) = 1$
when $(z,z^*)\in\Gamma_{m_k}$, thus $\|z\|_k \ge C^{-2}$.
Finally we know that every block subspace $Z$ of $Y$ contains
elements from $Q\Gamma_{m_k}$. If $z\in Z \cap Q\Gamma_{m_k}$,
we obtain for every $j\ne k$
$$ \|z\|_j \le 2 C^4 \alpha(2^{- \min(j,k)}) \|z\|_k.$$
This ends the proof of our remark.
\bigbreak

 We end with some comments and questions. It follows from our remark
that a uniformly convex Banach space contains an arbitrarily
distortable subspace, provided it contains an unconditional
basic sequence. According to a recent result of W.T.~Gowers [G],
the question of whether every uniformly convex Banach space
contains an arbitrarily distortable subspace
reduces now to the case of a uniformly convex
HI space (a Banach space $X$ is called HI --for {\it hereditarily
indecomposable}-- if the unit sphere of
every infinite dimensional subspace is
an asymptotic set in $X$); actually $X$ can also be 
assumed to be asymptotically-$\ell_p$;
is it possible for an HI space to be asymptotically-$\ell_p$?
\smallskip

 It is not so natural to try to find a distortable subspace rather than
trying to distort the whole space. If $X$ is a uniformly
convex Banach lattice, is it arbitrarily distortable? We did not
exactly prove this, even for a space with unconditional basis.
\smallskip

 Finally recall a problem raised by Milman-Tomczak: if $X$
does not contain $\ell_1$ or $c_0$, does it contain
an arbitrarily distortable subspace (or more: is it
arbitrarily distortable)? A first step in that direction
would be to elucidate the case of the Tsirelson space $T$
or of its dual $T^*$.
\bigbreak

\centerline{\bf References}
\medskip

\r{CS}{ P.G. Casazza, T.J. Shura, {\it Tsirelson's space},
Lecture Notes in Math. vol. 1363 (1989).}
\smallskip

\r{G}{ W.T.~Gowers, {\it A new dichotomy 
for Banach spaces}, preprint.}
\smallskip

\r{GM}{ W.T.~Gowers, B.~Maurey, {\it The unconditional basic
sequence problem}, Journal of AMS, to appear.}
\smallskip

\r{J}{ R.C.~James, {\it Uniformly non-square Banach spaces},
Ann. of Math. 80 (1964), 542--550.}
\smallskip

\r{K}{ J.L. Krivine, {\it Sous-espaces de dimension finie
des espaces de Banach r\'eticul\'es}, Ann. of Math. 104 (1976), 1--29.}
\smallskip

\r{LT1}{ J. Lindenstrauss, L. Tzafriri, {\it Classical Banach
spaces I: sequence spaces}, Springer Verlag, 1977.}
\smallskip

\r{LT2}{ J. Lindenstrauss, L. Tzafriri, {\it Classical Banach
spaces II: function spaces}, Springer Verlag, 1979.}
\smallskip

\r{M1}{ V.D.~Milman, {\it Spectrum of continuous bounded functions
on the unit sphere of a Banach space}, Funct. Anal. and Appl.
3 (1969), 67--79 (translated from Russian).}
\smallskip

\r{M2}{ V.D.~Milman, {\it The geometric theory of Banach spaces, part II:
Geometry of the unit sphere}, Uspekhi Math. Nauk 26 (1971), 73--149.
English translation in Russian Math. Surveys 26 (1971), 79--163.}
\smallskip

\r{MT}{ V.D.~Milman, N.~Tomczak-Jaegermann, {\it Asymptotic $\ell_p$
spaces and bounded distortions}, preprint.}
\smallskip

\r{OS}{ E.~Odell, T.~Schlumprecht,
{\it The distortion problem}, preprint.}
\smallskip

\r{S}{  T.~Schlumprecht,
{\it An arbitrarily distortable Banach space}, Israel J. Math. 76
(1991), 81--95.}
\smallskip

\r{T}{  B.S.~Tsirelson,
{\it Not every Banach space contains $\ell_p$ or $c_0$},
Funct. Anal. Appl. 8 (1974), 138--141 (translated from Russian).}
\smallskip

\bye